\documentclass[12pt, oneside, reqno]{amsart}

\usepackage{amsthm,amssymb,amsmath,amsfonts, tikz, graphicx, a4wide,subfigure,bm}
\usepackage{hyperref}
\usepackage{calligra}
\usepackage{color}
\usepackage{mathtools}
\usepackage{cite}
\usepackage{url}
\usepackage{comment}

\newtheorem{theorem}{Theorem}[section]
\newtheorem{proposition}[theorem]{Proposition}
\newtheorem{definition}[theorem]{Definition}

\newtheorem{lemma}[theorem]{Lemma}

\newtheorem{Theorem}{Theorem}

\newcommand{\ncm}{\newcommand}

\ncm{\be}{\begin{equation}}
\ncm{\ee}{\end{equation}}

\ncm{\ba}{\begin{array}}
\ncm{\ea}{\end{array}}

\ncm{\bea}{\begin{eqnarray}}
\ncm{\eea}{\end{eqnarray}}
\ncm{\bean}{\begin{eqnarray*}}
\ncm{\eean}{\end{eqnarray*}}

\ncm{\bt}{\begin{theorem}}
\ncm{\et}{\end{theorem}}

\ncm{\blm}{\begin{lemma}}
\ncm{\elm}{\end{lemma}}

\ncm{\bd}{\begin{definition}}
\ncm{\ed}{\end{definition}}

\ncm{\bp}{\begin{problem}}
\ncm{\ep}{\end{problem}}

\ncm{\lb}{\label}
\ncm{\lt}{\left}
\ncm{\rt}{\right}

\ncm{\f}{\frac}
\ncm{\RR}{\mathbb R}
\ncm{\NN}{\mathbb N}
\ncm{\mC}{\mathcal C}
\ncm{\sa}{\sigma}
\ncm{\ve}{\varepsilon}
\ncm{\rr}{\rho}
\ncm{\ta}{\theta}
\ncm{\al}{\alpha}

\let\eps=\varepsilon

\DeclareMathOperator{\diag}{diag}
\DeclareMathOperator{\id}{Id}

\begin{document}

\title[Differential Game]{On a linear differential game in the Hilbert space 
$\ell^2$}

\author[M.Ruziboev]{Marks Ruziboev}
\address[Marks Ruziboev]{Faculty of Mathematics, University of Vienna, Oskar-Morgenstern Platz 1, Austria}
\email{marks.ruziboev@univie.ac.at}

\author[K. Mamayusupov]{Khudoyor Mamayusupov}
\address[Khudoyor Mamayusupov]{ V.I.Romanovskiy Institute of Mathematics, Uzbekistan Academy of Sciences, 4, University street, Olmazor, Tashkent, 100174, Uzbekistan}
\email{mamayusupov@mathinst.uz}

\author[G.I. Ibragimov]{Gafurjan Ibragimov}
\address[Gafurjan Ibragimov]{University of Digital Economics and Agrotechnologies, Tashkent, 100022, Uzbekistan}
\email{ibragimov@upm.edu.my}

\author[A. Khaitmetov]{Adkham Khaitmetov}
\address[Adkham Khaitmetov]{Department of Higher Mathematics, Tashkent University of Information Technologies named after Muhammad Al-Khwarizmi, 100084, Tashkent, Uzbekistan}
\email{adhamhaitmetov@gmail.com}

\begin{abstract}
Two player pursuit evasion differential game and time optimal zero control problem in $\ell^2$ are considered. Optimal control for the corresponding zero control problem is found. A strategy for the pursuer that guarantees the solution for pursuit problem is constructed. 
\end{abstract}
\maketitle

\section{Introduction and statement of the problem}
Let for every $i\in \NN$, $ d_i\times d_i$ matrices $A_i$ be given with $2\le d_i\le d$, where $d\ge 2$ is a fixed integer.  We consider a differential game described by the following countable system of  differential equations 
\begin{equation}\label{eq:gsyst}
\dot x_i=A_i x_i+ u_{i}-v_{i},  x_i(0)=x_{i0}\in \RR^{d_i},  i=1,2,...,
\end{equation}
where $u_i$ is a control parameter of the pursuer, $v_i$ is a control parameter of the evader, both assumed to be  locally integrable functions with values in $\RR^{d_i}$ for all $i\in\NN$ satisfying certain constraints (see Definition \ref{def:constr}). For convenience we form column vector  $x=(x_1^\ast, x_2^\ast, \dots)^\ast$, where $x_i^\ast$ is the transpose of $x_i\in\RR^{d_i}$, for $i=1,2, \dots$\footnote{Actually, we are just concatenating vectors $x_1$, $x_2$, $\dots$ one below another to obtain an infinite vector. Below we adopt this point of view, which simplifies the notation considerably.}. Assume that $x_0=(x_{10}, x_{20}, \dots, )\in\ell^2$, i.e. $\|x\|^2=\sum_{i=0}^\infty\|x_{i0}\|^2<+\infty$, where $\|x_{i0}\|$ is the Euclidean norm of $x_{i0}\in\RR^{d_i}$. Let  $A=\text{diag}(A_1, A_2, ...)$ be an operator on $\RR^\infty$ whose action is defined as $Ax =\text{diag} ([A_1x_1]^\ast, [A_2x_2]^\ast, \dots)$ for $x\in\RR^\infty$. We would like to define $e^{tA}$ for suitable classes of matrices $A_i$. The problem here is that $A$ is not necessarily defined on $\ell^2$ or even on $\ell^\infty$ (i.e. $Ax$ is not necessarily in $\ell^\infty$ for $x\in \ell^\infty$ since $|A_ix_i|$ may go to infinity as $n\to \infty$). Therefore, we have to justify the existence of solutions to the above Cauchy problem for initial points in $\ell^2$.

We consider pursuit evasion differential game which consists of two separate problems as usual. The pursuit game can be completed in time $T>0$  provided there exists a  control function of the pursuer $u:\RR\to \ell^2$ such that for any control of $v:\RR\to \ell^2$ the solution of $x:\RR\to \ell^2$ of \eqref{eq:syst} for any $x_0$ satisfies $x(T)=0$. In this case $T$ is called the guaranteed pursuit time. Below we state precise conditions that are imposed on $u$, $v$.

A motivation for the setup comes from control problems for evolutionary PDEs, where using suitable decomposition of the control problem (see for example \cite{Cher95, Cher, Cher90, AzRuz, AzBakAkh}) $x_i$ would be a Fourier coefficients of an unknown function, while $u_i$ and $v_i$ would be that of  control parameters. Also, the setup is of independent interest as a controlled system in a Banach space (for works in this spirit see for example \cite{CurZwart,  AIMR}). Differential games for infinite dimensional systems are also well studied, for example, when the evolution  of the system is governed by parabolic equations pursuit evasion problems are considered in \cite{TM08, ST07, ST05}, where  the problem for the partial differential equations is reduced to an infinite system of  ordinary differential equations.  Pursuit evasion games with many players considered in \cite{Idham_IGI_Ask2016, Ibragimov-Al-Kuch2014, Ibragimov-ScAs2013, IFRP, TM08}.

For us the system \eqref{eq:gsyst} is a toy model of a system consisting of countably many point masses moving in $\RR^n$ with simple motion which are not interacting with each other. It is the first step in understanding the system of weakly interacting controllable particles in a more natural setting, e.g. for considering control problems for systems considered in  \cite{DF}.  

\section{Main results}\label{sec:sofprob}
As we pointed out earlier,  we have to justify the existence of solutions of the following Cauchy problem 
\begin{equation}\label{eq:syst}
\dot x_i=A_i x_i+ w_{i},  x_i(0)=x_{i0}\in \RR^{d_i},  i=1,2,...,
\end{equation}
with $w_i:\RR\to \RR^{d_i}$ locally integrable.
We look for  solutions of \eqref{eq:syst} from the space of continuous functions  $\mC([0,T]; \ell^2)$ for some $T>0$, such  that the coordinates $x_i(\cdot)$ of  $x:[0,T]\to \ell^2$ are almost everywhere differentiable.

\begin{definition}\label{def:constr}
We say that a family of matrices $\{A_i\}_{i\in\NN} $ is uniformly normalizable if there exists a family $\{P_i\}_{i\in\NN}$ of non-singular matrices and a constant $C\ge 1$ such that   $\|P_i\|\cdot\|P_i^{-1}\|\le C$ and  $P_iA_iP_i^{-1}$ is a matrix in the Jordan normal form for all $i\in \NN$. 
\end{definition}
Notice that there exist uniformly normalizable families of matrices, i.e., if all elements are already in Jordan normal form, then we may take $P_i=\text{Id}$ for all $i\in\NN$. On the other hand one can construct families, which aren't uniformly normalizable. 
In this work we will assume that the family of matrices in \eqref{eq:syst} and \eqref{eq:gsyst} are uniformly normalizable and  control parameters of the players satisfy the following constraint. 
\begin{definition}\label{def:constr}
Fix $\theta>0$ and let $B(\theta)$ be the set of all functions $w(\cdot)=\left(w_1(\cdot),w_2(\cdot),...\right),$ $w:[0,T] \to\ell^2$, with measurable coordinates $w_i(\cdot)\in \RR^{d_i}$, $0 \le t \le T$, $i=1,2,...$, that satisfy the constraint
\begin{equation}\label{eq:constr}
\sum\limits_{i=1}^\infty\int\limits_0^T\|w_{i}(s)\|^2ds \le \theta^2.
\end{equation}
$B(\theta)$ is called the set of admissible control functions.
\end{definition}
We have the following
\begin{Theorem}
\label{thm:ext}  Let $\{A_i\}$ be a family of uniformly normalizable matrices. 
If the real parts of eigenvalues of matrices $A_1, A_2,..$ are negative, then for any $w\in B(\theta)$, $\theta>0$ system \eqref{eq:syst} has a unique solution for any $z_0\in \ell_2$. Moreover, the corresponding components of the solution $x(t)=(x_1(t), x_2(t), \dots)$ are given by
\begin{equation}\label{eq:solution}
x_i(t)= e^{tA_i}x_{i0}+\int\limits_{0}^{t} e^{(t-s)A_i}w_i(s)ds, \quad i\in\NN.
\end{equation}
\end{Theorem}

\begin{definition}\label{def:0cont}
System  \eqref{eq:syst} is called \textbf{globally asymptotically stable} if $\lim_{t\to +\infty} x(t)=0$  for a solution $x(t)$ of  \eqref{eq:syst} with any initial condition $x_0\in\ell^2$ and $w_i\equiv 0$ for all $i\in \NN$. Further, system  \eqref{eq:syst} is \textbf{null-controllable} from $x_0\in\ell^2$ if there exists an admissible control $ u\in B(\theta)$  and $T=T(u)\in\RR_+$ such that the solution of  \eqref{eq:syst} starting from $x_0$ satisfies  $x(T)=0$. We say that system  \eqref{eq:syst} is \textit{null-controllable in large} if it is null-controllable from any $x_0\in\ell^2$. Also, $\inf_{u\in B(\theta)}T(u)$ is called \textbf{optimal time  of translation} and $u\in B(\theta)$ realizing the minimum is called \textbf{time optimal control}.
\end{definition}

\begin{Theorem}\label{thm:stability}
Under the assumptions of Theorem \ref{thm:ext} system \eqref{eq:syst} is globally asymptotically stable and  null controllable in large. Time optimal control exists and can be constructed explicitly. 
\end{Theorem}

Notice that the explicit form of the time optimal control requires some preliminary contraction and it is given in Section \ref{sec:proof_of_control}.

Further, we consider a pursuit-evasion differential game \eqref{eq:gsyst}. Fix $\rho, \sigma>0$. A function $u(\cdot)\in B(\rho)$ ($v(\cdot)\in B(\rho)$) is called an admissible control of the pursuer (evader). 

\begin{definition}\label{def:strategy}
A function $u:[0,T]\times \ell^2\to \ell^2$ with coordinates $u_k(t)=v_k(t)+\omega_k(t)$, $\omega\in B(\rho-\sigma)$, which is an admissible control of the pursuer for every $v\in\ell^2$ is called a strategy of the pursuer. 
\end{definition}

\begin{Theorem}\label{thm:game}
Suppose that $\rho>\sigma$ and the assumptions of Theorem \ref{thm:ext} are satisfied. Then, for any admissible control of the evader $v$  there exists a strategy of the pursuer $u$ and $\vartheta_1>0$ such that the solution of \eqref{eq:gsyst} satisfies $z(\tau)=0$ for some $0\le\tau\le \vartheta_1$, i.e., the game \eqref{eq:gsyst} can be completed within time $\vartheta_1$. 
\end{Theorem}
\section{Existence and uniqueness}
Notice that if we define $x(t)=(x_1(t), x_2(t), \dots)$ by setting every component $x_i$ as in  \eqref{eq:solution}, then $x(t)$ satisfies the equation and initial conditions in \eqref{eq:syst}.  This also implies uniqueness of a solution.   Thus it is sufficient to prove that $x(\cdot) \in \mathcal C([0, T], \ell^2)$  for any $T>0$.  Now we will show that $x(t)\in\ell^2$ for all $t\ge 0$. 
\subsection{Estimate for $\|e^{tA}\|$} 
Since $A=\diag\{A_1, A_2, \dots\}$ we have $e^{tA}=\diag\{e^{tA_1}, e^{tA_2},\dots\}$. Recall that for every $i$ there exists a non-singular transformation $P_i:\RR^{d_i}\to \RR^{d_i}$ such that
$ A_i=P_iJ_iP_i^{-1}$, where $J_i$ is the Jordan normal form of $A_i$. 

Thus 
\[
\|e^{tA_i}\|\le \|P_i\|\cdot \|e^{tJ_i}\|\cdot\|P_i^{-1}\|\le C\|e^{tJ_i}\|.
\]
By the assumption all eigenvalues $\lambda_{i1},\dots, \lambda_{id_i}$ of matrices have negative real part, letting $2\alpha_i=-\max_{1\le j\le d_i}Re(\lambda_{jd_i})>0$, we can find a polynomial $Q_i(t)$ of degree at most $d_i\le d$ (see \cite[\S 13]{Dech}) such that 
\begin{equation}\label{eq:|etai|}
\|e^{tA_i}\|\le  C|Q_i(t)| e^{-2t\alpha_i}\le \bar Ce^{-t\alpha_i}.
\end{equation}
Thus, for any $x\in \ell^2$ and $t\in[0,+\infty)$  we have
\begin{equation}\label{eq:|eta|}
\|e^{tA}x\|^2=\sum_{i=1}^\infty \|e^{tA_i}x_i\|^2\le \sum_{i=1}^\infty \|e^{tA_i}\|^2\|x_i\|^2\le \bar C^2\|x\|^2.
\end{equation}
This implies that $e^{tA}:\ell^2\to \ell^2$ is a bounded linear operator for every $t\in[0, +\infty)$.
Also, it is standard to check that $e^{tA}$ is a semigroup, i.e., $e^{(t+s)A}=e^{tA}e^{sA}$. 
\subsection{Proof of theorem \ref{thm:ext}}
We start by showing that $x(t)\in \ell^2$ for all $x_0\in\ell^2$ and  for all $t\in [0, T]$. Indeed, 
\begin{equation}\label{eq:xtnorm}
\|x(t)\|^2\le \sum_{i=1}^\infty\|e^{tA_i}x_{i0}+\int\limits_{0}^{t} e^{(t-s)A_i}w_i(s)ds\|^2\le 2\|e^{tA}x_0\|^2+2\sum_{i=1}^\infty\|\int\limits_{0}^{t} e^{(t-s)A_i}w_i(s)ds\|^2. 
\end{equation}
Let us estimate the last term of the above inequality.   We have
\begin{equation}\label{eq:xtlt}
\|\int_{0}^{t} e^{(t-s)A_i}w_i(s)ds\|^2\le \bar C^2\left(\int_{0}^{T} \|w_i(s)\|ds\right)^2\le \bar C^2 T\left(\int_{0}^{T} \|w_i(s)\|^2ds\right)
\end{equation}
where in the last step we have used the Cauchy-Schwartz inequality for $1$ against $\|w_i(s)\|$.
First substituting \eqref{eq:xtlt} into \eqref{eq:xtnorm} and then using \eqref{eq:|eta|}, \eqref{eq:|etai|} and constraint \eqref{eq:constr}  we  obtain 
\[
\|x(t)\|^2\le 2\bar C^2\|x_0\|^2+2\bar C^2T\theta^2,
\]
which proves the claim. 

We are now ready to prove that $x(t)=(x_1(t), x_2(t), \dots)\in \mC([0,T], \ell^2)$ for any $T>0$. 
Since $\|x(t)\|^2$ is bounded by a constant independent of $t$, for any $\eps>0$ there exists $N=N(\eps, t, t_0)\in\NN$  such that 
\begin{equation}\label{eq:tail}
\sum_{i=N+1}^\infty\|x_i(t)-x_i(t_0)\|^2\le \frac\eps2.
\end{equation}
For any $t,t_0\in [0,T]$ with $t_0\le t$ we have 
\begin{align*}
&\sum_{i=1}^{N}\|x_i(t)-x_i(t_0)\|^2
\le \sum_{i=1}^{N}\|e^{tA_i}-e^{t_0A_i}\|^2\cdot\|x_{i0}\|^2+\\
&+\sum_{i=1}^N \left\|e^{tA_i}\int_{0}^te^{-sA_i}w_i(s)ds
-e^{t_0A_i}\int_{0}^{t_0}e^{-sA_i}w_i(s)ds\right\|^2=(I)+(II).
\end{align*}
We start by estimating $(I)$. Notice that 
\begin{equation}\label{eq:et-et0}
\|e^{tA_i}-e^{t_0A_i}\|\le \|P_i\|\cdot\|P_i^{-1}\|\cdot\|e^{tJ_i}-e^{t_0J_i}\| \le  \bar C\|e^{t_0J_i}\|\cdot |Q_i(t-t_0)| e^{-|t-t_0|\alpha_i}. 
\end{equation}
Recall that $Q_i(t-t_0)$ is a polynomial of degree at most $d$ with coefficients depending only on the dimension $d_i$ of $J_i$. Thus, we can find $\delta$ independent of $i$ such that 
\begin{equation}\label{eq:Qi}
\bar C^3 |Q_i(t-t_0)|\cdot \|x_0\|^2 <\frac{\varepsilon}{4}\text{ and }\bar C^3 |Q_i(t-t_0)|<\frac{\varepsilon}{8T}
\end{equation}
for all $t\in (t_0-\delta, t_0+\delta)$.
Thus by \eqref{eq:|etai|} and \eqref{eq:Qi}  we have the following estimate for $(I)$:
\begin{equation}\label{eq:head}
\begin{aligned} 
\sum_{i=0}^{N}\|e^{tA_i}-e^{t_0A_i}\|^2\cdot & \|x_0\|^2\le \sum_{i=0}^{N}\|e^{tA_i}\|^2\cdot\|e^{(t-t_0)A_i}-\id\|^2\cdot \|x_0\|^2\le\\
&\sum_{i=0}^{N}\bar C^3 |Q_i(t-t_0)| e^{-|t-t_0|\alpha_i}\cdot \|x_0\|^2 <\frac{\varepsilon}{4}.
\end{aligned}
\end{equation}
For $(II)$ we write 
\[
\sum_{i=1}^N \left\|(e^{tA_i}-e^{t_0A_i})\int_{0}^te^{-sA_i}w_i(s)ds 
-e^{t_0A_i}\left(\int_{t_0}^{t}e^{-sA_i}w_i(s)ds\right)\right\|^2.
\]
Thus for every $i$ every summand of the above sum is bounded by 
\[
\left(\int_{0}^t \|e^{(t-s)A_i}-e^{(t_0-s)A_i}\| \cdot \|w_i(s)\|ds 
+ \int_{t_0}^{t}\|e^{(t_0-s)A_i}\|\cdot\|w_i(s)\|ds\right)^2.
\]
Applying inequality \eqref{eq:et-et0} and \eqref{eq:Qi} to the first summand and  inequality \eqref{eq:|etai|} to the second we obtain 
\[
 \sum_{i=1}^N \left(\frac{\varepsilon}{8T}\int_{0}^T 1_{[0, t]}\cdot \|w_i(s)\|ds 
+ \int_{0}^{T}1_{[t_0, t]}\cdot\bar Ce^{-t\alpha_i}\cdot\|w_i(s)\|ds\right)^2, 
\]
which is bounded by 
\[
\sum_{i=1}^N \left(\int_{0}^T \left(\frac{\varepsilon}{8T} 1_{[0, t]}+1_{[t_0, t]}\cdot\bar C\right)\|w_i(s)\|ds\right)^2.
\]
Now, using Cauchy-Schwartz inequality we obtain 
\[
\left(\frac{\varepsilon}{8}  + |t_0-t| \cdot\bar C\right)^2 
\sum_{i=1}^N \int_{0}^T \|w_i(s)\|^2ds.
\]
Since $|t-t_0|<\delta$ choosing $\delta>0$ sufficiently small and using \eqref{eq:constr} we bound the latter expression by $\frac{\varepsilon}{4}$ for all $\eps <4\theta^2$. Therefore, we conclude $(II)<\frac{\varepsilon}{4}$. Combining this, estimate \eqref{eq:head} and \eqref{eq:tail} imply that we obtain that $x(\cdot)\in \mC([0,1], 1)$, hence finishes the proof. 
\section{Proof of theorem \ref{thm:stability}}\label{sec:proof_of_control} 
\subsection{Asymptotic Stability}
We will show that $\|x(t)\|\to 0 $ as $t\to \infty$. Since $x_0\in\ell^2$ for any $\eps$ there exists $N=N(\eps)$ such that $\sum_{i=N+1}^\infty\|x_{i0}\|^2<\frac{\eps}{2\bar C}$.
By \eqref{eq:solution} and \eqref{eq:|etai|} we have 
\[
\|x(t)\|^2=\sum_{i=1}^\infty\|x_i(t)\|^2=\sum_{i=1}^\infty\|e^{tA_i}x_{i0}\|^2\le \sum_{i=1}^{N}\bar Ce^{-t\alpha_i}\|x_{i0}\|^2+\bar C\sum_{i=N+1}^\infty\|x_{i0}\|^2.
\]
Letting $\alpha_{\min}=\min_{1\le i\le N}\alpha_{i}>0$, from the above inequality we obtain
\[
\|x(t)\|^2\le \sum_{i=1}^{N}\bar Ce^{-t\alpha_i}\|x_{i0}\|^2+\bar C\sum_{i=N}^\infty\|x_{i0}\|^2\le \bar Ce^{-t\alpha_{\min}} \|x_0\|+\frac{\eps}{2}.
\]
There exists $t_\eps$ such that $\bar Ce^{-t\alpha_{\min}} \|x_0\|\le \frac{\eps}{2}$ for all $t<t_\eps$. 
This finishes the proof. 

Notice that if $\alpha_{\inf}=\inf_{i\ge 1}\alpha_{i}>0$, then the system is exponentially stable. Since in this case we don't have to cut at $N$, and can write 
\[
\|x(t)\|^2\le \sum_{i=1}^{\infty}\bar Ce^{-t\alpha_i}\|x_{i0}\|^2\le \bar Ce^{-t\alpha_{\inf}}\sum_{i=1}^{\infty}\|x_{i0}\|^2\le \bar Ce^{-t\alpha_{\inf}}\|x_{0}\|^2.
\]

\subsection{Gramians}
In this subsection we prove null controllability  of \eqref{eq:gsyst} and hence the proof of Theorem \ref{thm:stability}.
Our approach relies on Gramian operators and an observability inequalities. Set 
\[
W(\tau)=\int_0^\tau e^{-sA}\cdot e^{-sA^\ast}ds, \quad \tau\in \RR
\]
where $A^\ast$ is the adjoint of $A$ in $\ell^2$. The definition implies
\[
W(\tau)=\diag\{W_1(\tau), W_2(\tau), W_3(\tau)\dots \}, \text{ with }W_i(\tau)=\int_0^\tau e^{-sA_i}\cdot e^{-sA_i^\ast}ds, i\ge 1.
\]
Since $A_i$ is a finite matrix for every $i$ we have that $W_i(\tau)$ is a positive definite, symmetric and invertible operator\footnote{Notice that $W(\tau)$ is not necessarily bounded operator for fixed $\tau\in\RR$.} for every $i$ and $\tau\in\RR$, i.e., $W_i^{-1}(\tau)$ exists and bounded.  Define 
\[
W^{-1}(\tau)=\diag\{W_1^{-1}(\tau), W_2^{-1}(\tau), W_3^{-1}(\tau)\dots \}
\]
It is clear that $W_1^{-1}(\tau)$ is inverse to $W_1(\tau)$ for each $\tau\in \RR$. We will show that $W^{-1}(\tau):\ell^2\to \ell^2$ is a bounded linear operator.
For every $s\in\RR$, $i\in\NN$ and $x_i\in \RR^{d_i}$ we have 
\[
\langle W_i(\tau)x_i, x_i\rangle =\int_0^\tau \|e^{-sA_i^\ast}x_i\|^2ds\ge \int_0^\tau (m(P_i)e^{\beta_i s}m(P_i^{-1})\|x_i\|)^2ds=\int_0^\tau\frac{e^{2\beta_i s}\|x_i\|^2}{\|P_i\|^2\cdot \|P_i^{-1}\|^2}ds.
\]
where $m(P_i)$ is the minimum seminorm of $P_i$ and $\beta_i=-\min_{1\le j\le d_i}Re(\lambda_j)$. Since the eigenvalues of $A_i$ assumed to have strictly negative real part bounded away from zero, we have $\beta=\inf_i\beta_i>0$.  Therefore we have 
\[
\|W_i(\tau)\| \ge \frac{\langle W_i(\tau)x_i, x_i\rangle}{\|x\|^2}\ge \frac{1}{C^2}\int_0^\tau  e^{2\beta_is}ds\ge \frac{1}{C^2\beta_i} \left(e^{2\beta_i\tau}-1\right),
\]
which implies 
\begin{equation}\label{eq:wi^-1}
\|W_i^{-1}(\tau)\|\le 2C^2\beta_i\left(e^{2\beta_i\tau}-1\right)^{-1}\le C^2/\tau.
\end{equation}
Further, for any $x=(x_1, x_2, \dots)\in \ell^2$ with $\|x\|=1$  we have 
\begin{equation}\label{eq:w^-1norm}
\|W^{-1}(\tau)x\|^2= \sum_{i=1}^\infty\|W_i^{-1}(\tau)x_i\|^2\le \sum_{i=1}^\infty\|W_i^{-1}(\tau)\|^2\cdot\|x_i\|^2\le  C^2/\tau.
\end{equation}

\subsection{Null controllability in large}
Below we assume that $\theta>0$ and the set of admissible control is defined as in Section \ref{sec:sofprob}. 
Recall that $x(t)=e^{tA}x_0 + e^{tA}\int_{0}^te^{-sA}w(s)ds$  is the unique solution of  system \eqref{eq:syst} with an initial state $x(0)=x_0$. It is standard to check that the function 
\begin{equation}\label{eq:control}
    u^0(t)=-e^{-tA^\ast}\cdot W^{-1}(\tau)x_0 \quad\text{for every}\quad x_0\in \ell^2, \tau\in\RR^+
\end{equation}
solves the control problem if it is admissible, i.e., $\int_0^\tau e^{-sA}u^0(s)ds=-x_0$ for every fixed $\tau\in \RR^+$. Indeed, by \eqref{eq:control} we have
\begin{equation}\label{eq:int=y_0}
-\int_0^\tau e^{-tA_i}u^0dt=
\int_0^\tau e^{-tA_i}e^{-tA_i^\ast}dt\cdot W_i^{-1}(\tau)x_{i0}= x_{i0}, \text{ for all } i\in\NN.
\end{equation}

Therefore it remains to show that $u^0$ is admissible,  i.e.  there exists $\tau>0$ such that  $\|u^0\|^2=\sum\limits_{i=1}^\infty\int\limits_0^\tau\|u_{i}^0(s)\|^2ds\le \theta^2$, $u_i^0(s)\in\RR^{d_i}$.

By definition of $W(\tau)$ and Chauchy-Schwarz inequality  we have
\begin{equation}\label{eq:norm}
\begin{aligned}
\int_0^\tau\|u^0(t)\|^2dt&=\int_0^\tau \|e^{-tA^\ast}W^{-1}(\tau)u^0\|^2dt\\
&= \int_0^\tau \left\langle e^{-tA}\cdot e^{-tA^\ast}W^{-1}(\tau)x_0, W^{-1}(\tau)x_0\right\rangle dt\\
&=\langle x_0, W^{-1}(\tau)x_0\rangle\le \|x_0\|^2\cdot\|W^{-1}(\tau)\|.
\end{aligned}
\end{equation}
This together with inequality \eqref{eq:w^-1norm} and \eqref{eq:control}  implies that $u^0$ is admissible if 
\begin{equation}\label{eq:time}
C\|x_0\|^2/\sqrt\tau\le \theta^2.
\end{equation}
This finishes the proof, since the left hand side of \eqref{eq:time} decays as $\tau$ grows. 
\subsection{Time optimal control} 
Equation \eqref{eq:w^-1norm} shows that $\langle x_0, W^{-1}(\tau)x_0\rangle$ is decreasing as $\tau$ for every $x_0\in\ell^2$. Thus, for every $x_0\in\ell^2$ there exists a unique $\vartheta\in\RR_+$ such that 
\begin{equation}\label{eq:vartheta}
\langle x_0, W^{-1}(\tau)x_0\rangle>\theta^2, \text{ for } \tau>\vartheta, \text{ and }  \langle x_0, W^{-1}(\vartheta)x_0\rangle=\theta^2.
\end{equation}
We claim that $\vartheta$ is the optimal time.  We use the following result from \cite{LiM}.
\begin{proposition}
Let $B(t)$, $t\in[0,\vartheta_0]$ be a  continuous matrix-function of the order $d$ with a determinant not identically $0$ on $[0, \vartheta_0]$. Then among the measurable functions $w:[0, \vartheta_0]\to \RR^d$, satisfying the condition $\int_0^{\vartheta_0}B(s)w(s)ds=w_0\in\RR^d$ the function defined  almost everywhere on $[0,\vartheta]$ by the formula
\(
w(s)=B^\ast F^{-1}(\vartheta_0)x_0,\) \( F(\vartheta_0)=\int_0^{\vartheta_0}B(s)B^\ast(s)ds 
\)
gives a minimum to the functional $\int_0^{\vartheta_0}|w(s)|^2ds$.
\end{proposition}
Assume that there is an admissible control $u(\cdot)$ defined on $[0,\vartheta)$ such that $x(\tau)=0$ for some $\tau<\vartheta$. Be definition we have 
\[
e^{\tau A_i}x_{0i}+\int_0^\tau e^{(\tau-s)A_i}u_i(s)ds=0 \text{ for all } i\in\NN.
\]
Since $e^{(\tau-s)A_i}$ is continuous matrix function we can apply the above proposition for every $i\in\NN$ and conclude that the functional 
\(
J(u)=\int_{0}^\tau\sum_{k=1}^\infty\|u_i(s)\|^2ds 
\)
is minimized by $u^0$ defined in \eqref{eq:control}. Thus we have 
\[
J(u)\ge J(u^0)=\int_{0}^\tau\sum_{k=1}^\infty\|u_i^0(s)\|ds 
=\langle x_0, W^{-1}(\tau)x_0\rangle>\langle x_0, W^{-1}(\vartheta)x_0\rangle=\theta^2.
\]
This shows that $u(\cdot)$ is not admissible. This contradiction implies that $\vartheta$ is the optimal time of translation to the origin and $u^0(t)=-e^{-tA^\ast}\cdot W^{-1}(\vartheta)x_0$ is the time optimal control.
\section{Differential game problem: Proof of theorem \ref{thm:game}}
We now consider the game problem \eqref{eq:gsyst}. Recall that the equation 
\[
\langle x_0, W^{-1}(\tau)x_0\rangle =(\rho-\theta)^2
\]
has a unique solution $\vartheta_1$. Fix $T>\vartheta$. 
We define 
\begin{equation}\label{eq:pursuitcont}
u(t,v)=v-e^{-tA^\ast}\cdot W^{-1}(\vartheta_1)x_0 
\end{equation}
Let $v(\cdot)$ be any admissible control of the evader. We show that \eqref{eq:pursuitcont} is admissible.
\[
\|u(t,v)\|=\|v\|+\|e^{-tA^\ast} W^{-1}(\vartheta_1)x_0\|\le \sigma+ \langle x_0, W^{-1}(\vartheta_1)x_0\rangle^{1/2}=\rho.
\]
Also, it is easy to show that $x(\vartheta_1)=0$. This completes the proof. 
\section{Conclusion}
In this paper we studied infinite controllable system consisting of independent finite dimensional blocks. We solved optimal zero control problem  and constructed guaranteed strategy for pursuer to complete the pursuit game. We use Gramians in order construct optimal control.  It would be more desirable to consider more general equation than \eqref{eq:gsyst},
But we left this for further investigation. Since, our results don't generalize to this setting, and also one needs to find an analogue of Kalmann condition on controllability. 

We proved that the pursuit game can be completed if $\rho<\sigma$. Since in our setting the system is globally asymptotically stable and we expect that it is possible to complete the pursuit game for any $\rho, \sigma>0$. However, it turned out to be a challenging problem to define the strategy for $0<\rho<\sigma$. 
Also, we didn't attempt here evasion problem. We think that in the interval $(0,\vartheta_1)$ evasion is possible. However we leave this for future work. 

\section*{Acknowledgement } The research of M. Ruziboev (MR) is supported by the Austrian Science Fund (FWF): M2816 Meitner Grant. 

\section*{Data availability statement} We didn't generate and analysis any data in this work. 

\bibliographystyle{plain}

\end{document}